\documentstyle[12pt]{article}
\newcommand{\Res}{\mathop{\rm Res}\nolimits}
\newcommand{\sgn}{\mathop{\rm sgn}\nolimits}
\newcommand{\tr}{\mathop{\rm tr}\nolimits}
\newcommand{\diag}{\mathop{\rm diag}\nolimits}

\begin{document}

\title{On similar matrices over the dual numbers}
\author{I.M. Trishin}
\date{ }

\maketitle

Abstract. Matrices
over the dual numbers are considered.
We propose an approach to classify these matrices up to  similarity.
Some preliminary results on the realization of this approach
are obtained. In particular, we produce
explicitly canonical matrices
of orders 2 and 3.

\setcounter{section}{0}
\section{Introduction}

Let $K$ be an algebraically closed field of characteristic
zero and ${\cal D}=K[\zeta]$ is the algebra of the dual numbers over $K$,
where $\zeta^2=0$.
The main problem we consider in this paper
is to classify  elements of the full matrix algebra
$M_n({\cal D})$ up to similarity.
To be more precise, our goal is to determine a set of
{\it canonical matrices} such that each class of similar matrices
contains exactly one canonical matrix.

In Section 2, first, we reduce the problem to the matrices
$A_0+A_1\zeta$, where the set of  eigenvalues of
$A_0\in M_n(K)$ contains a single element
(see \linebreak
Theorem 2.2; note that this result
generalizes Berezin's theorem
(see, for example, Section 4 in  \cite{Red})).
Moreover, it can be assumed that $A_0=J_{\nu}$,
where $J_{\nu}$ is a block diagonal matrix such that its
diagonal blocks are Jordan blocks with the zero
eigenvalues.

Let $\mu=(\mu_1,\ldots,\mu_l)$  be a sequence of positive
integers and $\mu_1+\cdots+\mu_l=m$. Suppose the block partition
of matrices of $M_m(K)$ is defined by the sequence $\mu$.
We say that matrices $C,\tilde C\in GL_m(K)$ are
{\it $\mu$-mutual} if $C$ is lower triangular,
$\tilde C$ is upper triangular and  respective
diagonal blocks  of $C$ and $\tilde C$  are equal.
Matrices $A,B\in M_m(K)$ are {\it $\mu$-similar}
if there exist \linebreak
$\mu$-mutual matrices $C,\tilde C$
such that $B=CA\tilde C^{-1}$.
Concluding the process of reduction, in Section 2,
in the context of the appropriate theory,
we show that the
solution of the main problem can be realized in two steps.
The first one is
to classify  matrices of $M_m(K)$
up to  $\mu$-similarity.

Also, in Section 2
we prove that,
by analogy with the classic case,
if a \linebreak
class of  similar
${\cal D}$-matrices contains a diagonal matrix,
then this matrix is \linebreak
determined uniquely up to permutation
of  diagonal elements \linebreak
(Theorem 2.1).

In Section 3 we classify the sets of  $\mu$-similar
matrices  for $\mu=(1,1,\ldots,1)$ (Theorem 3.1).

In Section 4 we consider the matrices that are similar
to a matrix of the form $J_{n,\alpha}+A_1\zeta$,
where $J_{n,\alpha}$ is the Jordan block with the eigenvalue
$\alpha\in K$. For an arbitrary class of matrices
of this kind, we choose a canonical matrix
in a specific way (Theorem 4.1).

Also in this Section using the possibilities of the proposed approach we obtain classification of elements of $M_n({\cal D})$ up to similarity for $n=2,3$ (see Examples 4.1 and 4.2).

The author would like to express the gratitude to P.A. Saponov and D.I. Gurevich for their friendly support.

\section{The reduction of the main problem}
\setcounter{equation}{0}


Along with the classic case
matrices $A,B\in M_n({\cal D})$ are called {\it similar} if there
exists an invertible matrix $C\in M_n( {\cal D})$ such that
$B=CAC^{-1}$. Then we write $A\sim B$. If $A,B\in M_n(K)$
are similar in the classic sense then we say these are
{\it $K$-similar} and write $A \stackrel {K}{\sim} B$.

First consider the result that is of interest itself.
It is well known that if a matrix $A\in M_n(K)$
is $K$-similar to a diagonal matrix $A^{\prime}\in M_n(K)$,
then the set of  the diagonal elements of $A^{\prime}$
is determined uniquely by $A$. This set is
made up of the roots
of the characteristic polynomial $\chi(t)$ of $A$.
On the other hand, in our case a polynomial
may have not a unique expansion into
coprime factors.
Hence, in particular, distinct diagonal matrices
may have the same characteristic polynomial.
For example, for all diagonal matrices of the form
$
\left(
\begin{array}{cc}
a\zeta&0\\
0&-a\zeta
\end{array}
\right) ,
$
where $a\in K$, we have
$\chi (t)=t^2$.
Nevertheless,
suppose $A=\diag (a_{11}, \ldots,a_{nn})$,
$A^{\prime}=\diag(a_{11}^{\prime},\ldots,a_{nn}^{\prime})$,
where $a_{ii},a_{ii}^{\prime}\in {\cal D}$, $i=1,2,\ldots,n.$
Then we have

{\bf Theorem 2.1.} {\em If $A$ and $A^{\prime}$ are similar,
then
$$
\{ a_{11},\ldots,a_{nn} \}=
\{ a_{11}^{\prime},\ldots,a_{nn}^{\prime} \}.
$$
}

{\sc Proof.}  Write the matrix $A$ in the form $A=A_0+
A_1\zeta$, where
\begin{equation}
A_0=
\diag(
\underbrace{\alpha_1,\ldots,\alpha_1}_{\vartheta_1},
\underbrace{\alpha_2,\ldots,\alpha_2}_{\vartheta_2},
\ldots,
\underbrace{\alpha_l,\ldots,\alpha_l}_{\vartheta_l}
),
\label{ba1}
\end{equation}
$\alpha_1,\ldots,\alpha_l$ are pairwise distinct elements of $K$;
$\vartheta_1,\ldots,\vartheta_l$ are positive integers such that
$\vartheta_1\ge \vartheta_2\ge\ldots\ge \vartheta_l$,
$\vartheta_1+\cdots+\vartheta_l=n$
(in other words, $\vartheta=(\vartheta_1,\ldots,\vartheta_l)$
is the partition of $n$); $A_1\in M_n(K)$.

Since $A,A_0$ are diagonal, $A_1$
is also diagonal. In the same way, $A^{\prime}=
A_0^{\prime}+A_1^{\prime}\zeta,$ where
$A_0^{\prime}, A_1^{\prime} \in M_n(K)$ are diagonal.
By assumption,

\begin{equation}
A^{\prime}=GAG^{-1},
\label{ba2}
\end{equation}
where $G\in GL_n({\cal D})$. Represent the matrix $G$ in the form
$G=(E_n+C_1\zeta)B$, where $B\in GL_n(K)$,
$C_1\in M_n(K)$, $E_n$ is the identity matrix of order $n$.
Then we have $G^{-1}=B^{-1}(E_n-C_1\zeta)$.
Now from (\ref{ba2}) it follows that
$$
A_0^{\prime}=BA_0B^{-1},
$$
$$
A_1^{\prime}=BA_1B^{-1}+C_1A_0^{\prime}
-A_0^{\prime}C_1.
$$
Without loss of generality it can be assumed that
$A_0^{\prime}=A_0$, that is,
\begin{equation}
BA_0B^{-1}=A_0.
\label{ba3}
\end{equation}
Then we obtain
\begin{equation}
A_1^{\prime}=BA_1B^{-1}+C_1A_0-A_0C_1.
\label{ba4}
\end{equation}

Suppose for any matrix we consider the block partition is
defined by
the partition $\vartheta$. It follows from
(\ref{ba1}) and (\ref{ba3}) that $B$ is block diagonal.
Since $A_1$ is diagonal,we see that $BA_1B^{-1}$ is also block
diagonal. Recall that $A_1^{\prime}$ is diagonal.
Then equation (\ref{ba4}) implies that the matrix
$C_1A_0-A_0C_1$ is block diagonal.
Hence, from (\ref{ba1}) it follows that $C_1$ is
block diagonal. Therefore we have
\begin{equation}
C_1A_0-A_0C_1=0.
\label{ba5}
\end{equation}
Since $B$ is block diagonal,
from (\ref{ba4}) and (\ref{ba5}) it follows that
the diagonal matrices $A_1$ and $A_1^{\prime}$ differ
from one another by a permutation of diagonal elements
within the blocks determined by the partition $\vartheta$.
With (\ref{ba1}) we see that the theorem is proved. 

Let
\begin{equation}
A=A_0+\zeta A_1
\label{b0}
\end{equation}
be a matrix of
$M_n({\cal D})$, where $A_0,A_1\in M_n(K).$
By $\Lambda(A_0)$ denote the set of all eigenvalues of $A_0$.
We have the following important result:

{\bf{Theorem 2.2.}} {\em{Suppose $A_0=A_0^{\prime}\dot +
A_0^{\prime\prime}\equiv \left(
\begin{array}{cc}
A_0^{\prime}&0\\
0&A_0^{\prime\prime}
\end{array}
\right)$
and \\
$\Lambda(A_0^{\prime})\cap \Lambda(A_0^{\prime\prime})=
\emptyset.$
Then there exists a matrix $B\in M_n({\cal D})$ such that

i) $B$ is block diagonal, i.e., $B=(B_0^{\prime}+
B_1^{\prime}\zeta)\dot + (B_0^{\prime\prime}+
B_1^{\prime\prime}\zeta)$, where  \\
$B_i^{\prime},B_i^{\prime\prime}\in M_n(K)$ for $i=0,1$;

ii) $B \sim A$;

iii) $A_0^{\prime} \stackrel{K}{\sim} B_0^{\prime}$ and
$A_0^{\prime\prime} \stackrel{K}{\sim} B_0^{\prime\prime}$.

}}

{\sc{Proof.}} Let $f=f_0+f_1\zeta,\: f_0^{\prime},\:
f_0^{\prime\prime}$ be the characteristic polynomials
of the matrices $A, A_0^{\prime}, A_0^{\prime\prime}$
respectively, where $f_0, f_1, f_0^{\prime},
f_0^{\prime\prime}\in K[t]$.
Since $A_0=A_0^{\prime}\dot + A_0^{\prime\prime}$,
we have $f_0=f_0^{\prime}f_0^{\prime\prime}$.
Since the polynomials $f_0^{\prime},f_0^{\prime\prime}$
are coprime,  there exist polynomials
$g^{\prime},g^{\prime\prime}\in K[t]$ such that
$$
f_0^{\prime}g^{\prime}+f_0^{\prime\prime}g^{\prime\prime}=1
$$
(see \cite{Lang}).
Now it can easily be checked that
$$
f=f^{\prime}f^{\prime\prime},
$$
where $f^{\prime}=f_0^{\prime}+
g^{\prime\prime}f_1\zeta,\:
f^{\prime\prime}=f_0^{\prime\prime}+
g^{\prime}f_1\zeta.$

There exist polynomials $h^{\prime},h^{\prime\prime}
\in {\cal D}[t]$ such that
$$
h^{\prime}f^{\prime}+h^{\prime\prime}f^{\prime\prime}=
\Res(f^{\prime},f^{\prime\prime})
$$
(see \cite{Lang}). Because $\Res(f^{\prime},
f^{\prime\prime})\in {\cal D}$, we have
$$
h^{\prime}f^{\prime}+
h^{\prime\prime}f^{\prime\prime}=a_0+a_1\zeta,
$$
where $a_0,\:a_1\in K$. We see that
$\Res(f^{\prime},f^{\prime\prime}) \Bigr|_{\zeta=0}
=\Res(f_0^{\prime},f_0^{\prime\prime}).$
Whence, $a_0=\Res(f_0^{\prime},f_0^{\prime\prime}).$
Since $f_0^{\prime},f_0^{\prime\prime}$ are coprime,
$a_0\ne 0$ and the element $a_0+a_1\zeta$
is invertible.

Without loss of generality it can be assumed that $A$ is a matrix
of a ${\cal D}$-linear transformation ${\cal A}$
of the ${\cal D}$-envelope $V_{\cal D}$ of a $K$-linear space $V$
with respect to some basis.

From our previous results (see Proposition 3.3 in \cite{Red})
it follows that
\begin{equation}
V_{\cal D}=\ker f^{\prime}({\cal A})
\oplus_{\cal D}
\ker f^{\prime\prime}({\cal A}).
\label{b1.0}
\end{equation}
Also, the ${\cal D}$-modules
$\ker f^{\prime}({\cal A})$ and
$\ker f^{\prime\prime}({\cal A})$
are free (see \cite{Red}).
By $B$ denote a matrix of ${\cal A}$ with respect
to some basis associated with  decomposition
(\ref{b1.0}). Then $B$ is block diagonal and
$A\sim B$.

The ${\cal D}$-linear transformation ${\cal A}$ of $V_{\cal D}$
determines uniquely the $K$-linear transformation
${\cal A}_0$ of the space $V$ for $\zeta=0$.
In this case, $A_0$ is the matrix of ${\cal A}_0$
with respect to some basis.
Suppose $V=V_1\oplus_KV_2$ is the decomposition of $V$
such that  the block diagonal form of
$A_0=A_0^{\prime}\dot +A_0^{\prime\prime}$
(see above) corresponds to this decomposition.
Let $f_0^{\prime}$,  $f_0^{\prime\prime}$ be as above.
Then we have

{\bf{Lemma 2.1.}}
$$
\ker f_0^{\prime}({\cal A}_0)=V_1,\;\;\;
\ker f_0^{\prime\prime}({\cal A}_0)=V_2.
$$

The proof of this lemma is left to the reader.

From Lemma 2.1 it follows that $A_0^{\prime}$,
$B_0^{\prime}$ ($A_0^{\prime\prime}$,
$B_0^{\prime\prime}$, respectively)
are the matrices of the same linear transformation
of the same linear space.
Hence we have
$A_0^{\prime}\stackrel{K}{\sim} B_0^{\prime}$
($A_0^{\prime\prime}\stackrel{K}{\sim} B_0^{\prime\prime}$).

This completes the proof of the theorem.  

Let $A$ be given by (\ref{b0}) and $|\Lambda(A_0)|>1.$
Then it is well known that $A_0$ is $K$-similar
to a block diagonal matrix
$B_0=B_0^{\prime}\dot + B_0^{\prime\prime}$
such that $\Lambda (B_0^{\prime})\cap \Lambda (B_0^{\prime\prime})
=\emptyset$. Therefore Theorem 2.2 yields the following result:

{\bf Corollary.} {\em If  $|\Lambda(A_0)|>1$, then the matrix
$A$ is similar to a block diagonal matrix.}

Thus, by induction, the main problem is reduced
to the case when  \linebreak
$|\Lambda(A_0)|=1$.
Since  diagonal matrices make up the center
of $M_n({\cal D})$, we can assume
that $\Lambda(A_0)
=\{ 0\}.$
Then the matrix $A$ is similar to a matrix $A^{\prime}$
of the form
$$
A^{\prime}=J_{\nu}+A_1^{\prime}\zeta,
$$
where

\begin{equation}
J_{\nu}=J_{\nu_1,0}\dot +J_{\nu_2,0}\dot +
\cdots \dot +J_{\nu_m,0},
\label{b1}
\end{equation}
$\nu=(\nu_1,\nu_2,\ldots,\nu_m)$ is a partition of $n$,
$J_{\nu_i,0}$ is the Jordan block
of order $\nu_i$ with the zero eigenvalues;
$A_1^{\prime}\in M_n(K)$.
Note that the matrix $J_{\nu}$ is determined uniquely
by the first term $A_0$ of $A$
(see (\ref{b0})).

Now, as for the main problem, it suffices
to choose a
canonical matrix in a certain way
for every set of
similar matrices of the form
$J_{\nu}+A_1\zeta$, where $A_1\in M_n(K)$.

Suppose
$$
D(J_{\nu} +A_1\zeta)D^{-1}=
J_{\nu} +A_1^{\prime}\zeta,
$$
where $D\in GL_n({\cal D});\;A_1,A_1^{\prime}\in M_n(K).$
Then we claim that the matrix $D$ can be represented
in the form
\begin{equation}
D=(E_n+C_1\zeta)B_0,
\label{b2}
\end{equation}
where
$C_1\in M_n(K)$, $B_0\in GL_n(K)$
and
$$
B_0J_{\nu} =J_{\nu} B_0.
$$
Indeed, suppose $D_0,D_1\in M_n(K)$ are matrices such that
$$
D=D_0+D_1\zeta.
$$
Then we see that
$$
D_0J_{\nu} =J_{\nu} D_0
$$
and
$$
D_0+D_1\zeta=(E_n+D_1D_0^{-1}\zeta)D_0.
$$
Therefore we have (\ref{b2}), where $C_1=D_1D_0^{-1},\;
B_0=D_0$.
Thus the following result is proved:

{\bf{Lemma 2.2.}} {\em{
If
$$
J_{\nu}+A_1\zeta \sim J_{\nu}+A_1^{\prime}\zeta,
$$
where $A_1,A_1^{\prime}\in M_n(K)$,
then there exists a matrix $D\in GL_n(K)$ of the form (\ref{b2})
such that
$$
J_{\nu}+A_1^{\prime}\zeta =D(J_{\nu}+A_1\zeta)D^{-1}. \; 
$$
}}

For an arbitrary positive integer $l$,
by $[l]$ denote the set $\{ 1,2,\ldots,l\}$.

Let $R,S$
be proper subsets of $[n]$;
$\bar R=[n]\setminus R$, $\bar S=[n]\setminus S$.
By $A(\bar R,\bar S)$ denote the submatrix
of  $A\in M_n(G_0)$ such that an element $a_{ij}$ of $A$ belongs to
$A(\bar R,\bar S)$ iff $i\in \bar R$ and $j\in \bar S$.
By definition, put $A_{R,S}=A(\bar R,\bar S)$.
Notice that $A_{R,S}$ is obtained from $A$ if we remove from $A$
all rows whose numbers belong to $R$
and all columns whose numbers belong to $S$.

Let the numbers of the non-zero rows and columns
of $J_{\nu}$ be $i_1,i_2,\ldots,i_r$ and $j_1,j_2,\ldots,j_r$
respectively; $P=\{i_1,\ldots,i_r\}$, $Q=\{j_1,\ldots,j_r\}$.

Suppose $A=J_{\nu}+A_1\zeta,\; C=E_n+C_1\zeta$,
where $A_1,\:C_1\in M_n(K)$; $A^{\prime}=CAC^{-1}$.

{\bf{Lemma 2.3.}}{\em{
$$
A^{\prime}_{P,Q}=A_{P,Q}.
$$
}}

{\sc{Proof.}} We have
$$
A^{\prime}=(E_n+C_1\zeta )A(E_n-C_1\zeta)=A+\hat A\zeta,
$$
where $\hat A=C_1J_{\nu}-J_{\nu}C_1.$
If an element $x_{ij}$ of the matrix $C_1J_{\nu}$
$(J_{\nu}C_1)$ is not equal to zero, then $j\in Q\;
(i\in P)$. Hence,
$$
\hat A_{P,Q}=0
$$
and we see that the lemma is proved.  

By definition, put
$$
\bar {\cal B}_{\nu}=
\{ B\in M_n(K)|BJ_{\nu}=J_{\nu}B\},
$$
$$
{\cal B}_{\nu}=
\{B\in \bar{\cal B}_{\nu}| \; |B|\ne 0\}.
$$

Evidently, ${\cal B}_{\nu}$ is a subgroup of $GL_n(K)$.

{\bf Remark.} Taking into account Lemmas 2.2 and 2.3
we see that now the first problem is to describe the sets of the form
$$
\{ (BA_1B^{-1})_{P,Q}| B\in {\cal B}_{\nu}\},
$$
where $A_1\in M_n(K)$. In the present section we prove
that if $B$ ranges over ${\cal B}_{\nu}$, then
$(BA_1B^{-1})_{P,Q}$ ranges over all matrices that are $\mu$-similar
to $(A_1)_{P,Q}$
(see Introduction and Theorem 2.4 below).

By $W$ denote the algebra $M_n(K)$ as a linear space.
Let $W_{\nu}$ be the subspace of $W$ such that
$A\in W_{\nu}$ iff $A_{P,Q}=0$.

For example, if $\nu=(2,1)$, then $W_{\nu}=
\left(
\begin{array}{ccc}
*&*&*\\
0&*&0\\
0&*&0
\end{array}
\right)=\langle e_{11},e_{12},e_{13},e_{22},e_{32}\rangle,$
where $e_{ij}$ are the matrix units.
It is clear that, generally, for any $\nu$ the  \linebreak
$K$-space $W_{\nu}$
has the basis consisting of some matrix units.

We see that $W$ is a ${\cal B}_{\nu}$-module with respect
to the adjoint action.

{\bf{Lemma 2.4.}}
$$
{\cal{B}}_{\nu}(W_{\nu})=W_{\nu}
$$

This lemma follows from

{\bf{Lemma 2.5.}} {\em{If $A\in W_{\nu}$ and $B\in \bar {\cal B}_{\nu}$,
then $BA\in W_{\nu}$ and $AB\in W_{\nu}$.}}

{\sc{Proof.}} First prove that $BA\in W_{\nu}$.
Assuming the converse, suppose $e_{ij}$ is
a matrix unit
such that

i) $e_{ij}$ belongs to $W_{\nu}$;

ii) $Be_{ij}=\alpha e_{kj}+\cdots$ , where $\alpha\in K,$
$\alpha\ne 0,\;e_{kj}\notin W_{\nu}$ and dots denote
a linear combination of matrix units that are distinct
from $e_{kj}$.
(Recall that $e_{kj}\notin W_{\nu}$ means that the $k$-th
row and the $j$-th column of the matrix $J_{\nu}$ are zero.)

Then
$$
B=\alpha e_{ki}+\cdots,
$$
where dots denote a linear combination of matrix units
that are distinct from $e_{ki}$.

Since the $k$-th row of $J_{\nu}$ is zero, the $k$-th row
of $J_{\nu}B$ is also zero. On the other hand, we claim that
the $k$-th row of $BJ_{\nu}$ is not zero. Indeed,
because $e_{ij}\in W_{\nu}$ and the $j$-th column of $J_{\nu}$
is zero, we see  that
the $i$-th row of $J_{\nu}$ is not zero. Hence there exists
a uniquely determined integer $l$ such that the element
$y_{pl}$ of the matrix $J_{\nu}$ is equal to zero for
$p\ne i$ and $y_{il}=1$.
Since $\alpha\ne 0$, we see that the $k$-th row of the matrix
$BJ_{\nu}$ is not zero.

Therefore, $J_{\nu}B\ne BJ_{\nu}$.
This contradiction concludes the proof.

Similarly, it can be proved that $AB\in W_{\nu}$.

This completes the proof of the lemma.  

By $\Phi$ denote the representation of the group ${\cal B}_{\nu}$
on the factor space $W/W_{\nu}$.

For an arbitrary proper subset $I$ of $[n]$, denote
by $A_{I}$ the matrix $A_{I,I}\equiv A(\bar I,\bar I)$,
 where $A\in M_n(K)$.

Let the sets $P$ and $Q$ be as above, $m=n-r$.
For an arbitrary $B\in {\cal B}_{\nu}$, denote by $\psi_B$ the mapping
of $M_m(K)$ to $M_m(K)$ such that
\begin{equation}
\psi_B(Z)=B_PZ(B_Q)^{-1}
\label{b7}
\end{equation}
for all $Z\in M_m(K)$.

Let the mapping $\Psi:{\cal B}_{\nu}\to End_K(M_m(K))$ be given by
$\Psi:B\mapsto \psi_B$.

{\bf{Lemma 2.6.}} {\em{The mapping $\Psi$ is a representation
of the group ${\cal B}_{\nu}$.}}

To prove this we need several lemmas.

Let $I$ be a subset of $[n]$, $A\in M_n(K)$. We say that the matrix
$A$ satisfies  condition $(I)$ if for all elements $a_{ij}$
of $A$ such that $i\in \bar I$ and $j\in I$ we have $a_{ij}=0$.

Let $A,B\in M_n(K).$

{\bf{Lemma 2.7.}} {\em{If the matrix $A$ satisfies   condition
$(I)$ or the matrix $B$ satisfies   condition $(\bar I)$,
then
$$
(AB)_I=A_IB_I.
$$
}}
{\sc{Proof.}} Let $c_{ij}$ be an element of the matrix $AB$.
Then $\displaystyle{c_{ij}=\sum_{t\in I\cup \bar I}a_{it}b_{tj}.}$
Suppose $i,j\in\bar I$, then from  condition $(I)$ for $A$
(or from  condition $(\bar I)$ for $B$) it follows that
$c_{ij}=\displaystyle{\sum_{t\in \bar I}a_{it}b_{tj}}$.
But the last formula
is just for the element of the matrix $A_IB_I$.

{\bf{Lemma 2.8.}} {\em{Any matrix $B\in \bar {\cal B}_{\nu}$
satisfies   conditions $(P)$ and $(\bar Q)$.}}

{\sc{Proof.}} Assume the converse, that is, suppose there
exists an element $b_{ij}$ of  $B$ such that
$i\in \bar P,\;j\in P,\;b_{ij}\ne 0$. Then there exists
a column of  $J_{\nu}$ such that its $j$-th element
is equal to $1$ and the rest elements are zero.
Hence the $i$-th row of the matrix $BJ_{\nu}$
is not zero. On the other hand, since the $i$-th row
of $J_{\nu}$ is zero, the $i$-th row of $J_{\nu}B$
is also zero. Thus we have $BJ_{\nu}\ne J_{\nu}B$.
This contradiction concludes the
proof of the first assertion of the
lemma.

Likewise, the second assertion can be easily proved.  

From Lemmas 2.7 and 2.8 it follows that for all
$B^{\prime},\:B^{\prime\prime}\in \bar{\cal B}_{\nu}$
we have
\begin{equation}
(B^{\prime}B^{\prime\prime})_P=
B^{\prime}_P B^{\prime\prime}_P,\;\;
(B^{\prime}B^{\prime\prime})_Q=
B^{\prime}_Q B^{\prime\prime}_Q.
\label{b3}
\end{equation}

{\sc{Proof of Lemma 2.6.}} In the mapping $B\mapsto B_P$
($B\mapsto B_Q$), if a row with a number $t$ is removed from
$B$ then the column with the same number is also removed.
Whence,
\begin{equation}
(E_n)_P=(E_n)_Q=E_m,
\label{b6}
\end{equation}
where $m=n-r$, and we have
$$
\Psi(E_n)=\psi_{E_n}=id.
$$
With (\ref{b3}) we obtain
$$
\psi_{B^{\prime}B^{\prime\prime}}(Z)=
(B^{\prime}B^{\prime\prime})_PZ\bigl( (B^{\prime}B^{\prime\prime})_Q
\bigr) ^{-1}=B_P^{\prime}B_P^{\prime\prime}Z
\bigl( B_Q^{\prime}B_Q^{\prime\prime} \bigr) ^{-1}
$$
$$
=B_P^{\prime}B_P^{\prime\prime}Z
\bigl( B_Q^{\prime\prime}\bigr)^{-1}
\bigl( B_Q^{\prime} \bigr) ^{-1}=
B_P^{\prime}\psi_{B^{\prime\prime}}(Z)
\bigl( B_Q^{\prime} \bigr)^{-1}=\psi_{B^{\prime}}
\bigl( \psi_{B^{\prime\prime}}(Z) \bigr),
$$
that is, we have

$$
\psi_{B^{\prime}B^{\prime\prime}}=
\psi_{B^{\prime}}\psi_{B^{\prime\prime}}. \;  
$$

{\bf{Proposition 2.1.}} {\em{The representations $\Phi$ and
$\Psi$ of the group ${\cal B}_{\nu}$ are equivalent
.}}

{\sc{Proof}}. Suppose matrices $A,C\in M_n(K)$ satisfy   conditions
$(R)$ and $(\bar S)$ respectively. Then we claim that
\begin{equation}
(ABC)_{R,S}=A_RB_{R,S}C_S,
\label{b4}
\end{equation}
where $B\in M_n(K)$. Indeed, let $d_{ij}$ be an element
of the matrix $ABC$. Then we have
$$
d_{ij}=\sum_{\textstyle{ k\in R\cup \bar R \atop l\in S\cup \bar S}}
a_{ik}b_{kl}c_{lj},
$$
where $a_{ik},b_{kl},c_{lj}$ are the elements of $A,B,C$
respectively. If $i\in \bar R$ and $j\in \bar S$, then
$$
d_{ij}=\sum_
{\textstyle{k\in\bar R\atop l\in \bar S}}
a_{ik}b_{kl}c_{lj}.
$$
But the same formula determines the element of the matrix
$A_RB_{R,S}C_{S}$. Hence equality (\ref{b4}) is proved.

Let $A\in M_n(K)$, $B\in {\cal B}_{\nu}$.

Using (\ref{b3}) and (\ref{b6}), we get
\begin{equation}
(B^{-1})_Q=B_Q^{-1}.
\label{b8}
\end{equation}

Then, by Lemma 2.8, (\ref{b4}) and (\ref{b8}), it follows that
\begin{equation}
(BAB^{-1})_{P,Q}=B_PA_{P,Q}B_Q^{-1}.
\label{b5}
\end{equation}

Let $\Theta$ be the mapping from $M_n(K)$ to $M_m(K)$
such that $\Theta(A)=A_{P,Q}$ for all $A\in M_n(K)$.
Obviously, the mapping $\Theta$ is linear.
Then with (\ref{b5}), (\ref{b7}) and (\ref{b8}) we obtain
$$
\Theta(BAB^{-1})=(BAB^{-1})_{P,Q}=B_PA_{P,Q}(B^{-1})_Q=
$$
$$
B_P\Theta(A)B_Q^{-1}=\psi_B(\Theta(A)).
$$
Thus the linear mapping $\Theta$ is a homomorphism
of the ${\cal B}_{\nu}$-modules.

Finally, by the definition of $W_{\nu}$,
we have
$$
\ker \Theta=W_{\nu}. \; 
$$

Let $\nu$ be as above. Write this partition in the form
\begin{equation}
\nu=(\underbrace{\alpha_1,\ldots,\alpha_1}_{s_1},
\underbrace{\alpha_2,\ldots,\alpha_2}_{s_2},\ldots,
\underbrace{\alpha_l,\ldots,\alpha_l}_{s_l}),
\label{b8.1}
\end{equation}
where $\alpha_i\ne \alpha_j$ for $i\ne j$,
$\alpha_i>\alpha_{i+1},\;s_i\ge 1,$ that is, $s_i$
is the multiplicity of the part $\alpha_i$ of $\nu$.
By $\hat\nu$ denote the sequence of the multiplicities
$s_i$ of $\nu$, that is,
$$
\hat\nu=(s_1,s_2,\ldots,s_l).
$$
For example, if $\nu=(5,5,4,4,4,2)$, then $ \hat\nu=(2,3,1)$.

In the following for any $B\in {\cal B}_{\nu}$
the block partition of the matrices $B_P,\:B_Q$
is defined by the
sequence $\hat\nu$.

Suppose $B\in {\cal B}_{\nu}$ is of the general form. Then we have

{\bf{Theorem 2.3.}} {\em{The matrix $B_P$ $(B_Q)$
is the invertible lower (respectively, upper) block triangular matrix
of the general form. Moreover

i)  respective diagonal blocks of $B_P$ and $B_Q$
are coincide;

ii) the  elements that belong to  off-diagonal
and non-zero blocks of $B_P$ and $B_Q$ are independent.}}

{\sc{Proof.}} We assume that
the block partition of the matrix $B$
is  defined by the partition $\nu=
(\nu_1,\nu_2,\ldots,\nu_m)$,
that is, $B=(B_{ij})$,
where $B_{ij}$ is a rectangular $\nu_i\times\nu_j$-matrix.

Recall from \cite{Gant} that a $k\times l$-matrix $G=(g_{ij})$ is
{\it regular} if the following conditions hold:

i)  $G$ is upper triangular;

ii) if $k>l$, then $g_{ij}=0$ for $i>l$;

iii) if $k<l$, then $g_{ij}=0$ for $j\le l-k$;

iv) $g_{ij}=g_{pq}$ if $i-j=p-q$.

A matrix $B\in M_n(K)$ belongs to ${\cal B}_{\nu}$ iff
every block  of $B$ is regular and $|B|\ne 0$.
(see, for example, \cite{Gant}).

By definition, put
$$
\sigma_k=\nu_1+\cdots+\nu_k,
$$
where $k=1,2,\ldots,m;\:\sigma_0=0$.

Consider a block $B_{ij}$ of $B$.

By definition,
$\sigma_t=\sigma_{t-1}+\nu_t\notin P$
and if $\nu_t>1$, then
$$
P\supseteq \{\sigma_{t-1}+1,\sigma_{t-1}+2,
\ldots,\sigma_{t-1}+\nu_t-1\},
$$
where $t\in [m]$.
Hence getting $B_P$ from $B$ we remove the first $(\nu_i-1)$
rows and the first $(\nu_j-1)$ columns of the block $B_{ij}$.
Thus  the only element of $B_{ij}$ is remained in $B_P$.
It is placed in the  right lower angle of $B_{ij}$. This element
we denote by $b_{ij}^{\prime}$.

By analogy, getting $B_Q$ from $B$ we see that
the only element of $B_{ij}$ is remained in $B_Q$.
It is placed in the left upper angle of $B_{ij}$.
This element we denote by $b_{ij}^{\prime\prime}$.

Let the block $B_{ij}$ be square. Since $B_{ij}$ is regular,
in particular,
we have, $b_{ij}^{\prime}=b_{ij}^{\prime\prime}$,
that is, the  respective elements of $B_P$ and $B_Q$ are equal.

Let the block $B_{ij}$ be rectangular. Suppose $i<j$.
Then, since $B_{ij}$ is regular, we have $b_{ij}^{\prime}=0$.
Therefore the matrix $B_P$ is lower block triangular.
In the same way, if $i>j$, then we have $b_{ij}^{\prime\prime}=0$.
Whence the matrix $B_Q$ is upper block triangular.

In the mapping $B\mapsto B_P$ $(B\mapsto B_Q)$ an aggregate
of all square blocks of the same order is taken to some
square diagonal block of $B_P$ $(B_Q)$.
Likewise, an aggregate of all rectangular blocks
of the same dimension is taken to the corresponding
block of $B_P$ ($B_Q$).

Now from Proposition 2.2 (see below) it follows that
the matrix $B_P$ ($B_Q$) is the invertible lower
(respectively, upper) block triangular matrix of the general form.

From the above, distinct elements of   off-diagonal
and non-zero blocks of $B_P$ and $B_Q$ are determined
by  distinct rectangular blocks of $B$.
Also, it follows from Proposition 2.2 that the determinant
$|B|$ depends only on the elements of the  square blocks of $B$.
Since $B$ is of the general form, we see that  elements
of distinct rectangular blocks of $B$ are independent.
This yields  statement ii) of the theorem.

{\bf Proposition 2.2.}
$$
|B|=0
\Longleftrightarrow
|B_P|=0.
$$

\vspace{0.5cm}

To prove this, we need several lemmas.

Suppose $\nu=(\underbrace{t,t,\ldots,t}_s)$,
where $st=n$.
As above, the block partition of $B\in \bar {\cal B}_{\nu}$
is defined by $\nu$,
that is, for any $i,j\in [s]$ the block $B_{ij}$ of $B$
is a square matrix of order $t$.
Recall that any block of $B$ is regular, i.e.,
$$
B_{ij}=b_{ij}^{(0)}E_t+
b_{ij}^{(1)}J_{t,0}+
b_{ij}^{(2)}J_{t,0}^2+
\cdots+
b_{ij}^{(t-1)}J_{t,0}^{t-1},
$$
where $b_{ij}^{(l)}\in K$; $i,j\in [s]$.

For example, if $n=4$, $t=2$, then
$$
B=
\left(
\begin{array}{cccc}
b_{11}^{(0)}&b_{11}^{(1)}&b_{12}^{(0)}&b_{12}^{(1)}\\
0&b_{11}^{(0)}&0&b_{12}^{(0)}\\
b_{21}^{(0)}&b_{21}^{(1)}&b_{22}^{(0)}&b_{22}^{(1)}\\
0&b_{21}^{(0)}&0&b_{22}^{(0)}
\end{array}
\right),
$$
$P=\{ 1,3\},$ $Q=\{ 2,4\}$.
Hence, in this case,
$$
B_P=B_Q=
\left(
\begin{array}{cc}
b_{11}^{(0)}&b_{12}^{(0)}\\
b_{21}^{(0)}&b_{22}^{(0)}
\end{array}
\right)
\equiv(b_{ij}^{(0)}).
$$
It can easily be checked that
$
|B|=|B_P|^2=|B_Q|^2.
$
By analogy, generally,
for $\nu=(\underbrace{t,t,\ldots,t}_s)$ we have

{\bf Lemma 2.9.}
$$
|B|=|B_P|^t=|B_Q|^t.
$$

{\sc Proof.} Let us prove that
\begin{equation}
|B|=|b_{ij}^{(0)}|^t.
\label{b9.0}
\end{equation}

 We proceed by induction on $t$. For $t=1$,
there is nothing to prove.

Let $t>1$. By  Laplace's theorem we have
\begin{equation}
|B|=
\sum_R
(-1)^{\sigma(R)+\sigma(S)}
|B(R,S)|
\cdot
|B(\bar R,\bar S)|,
\label{b9}
\end{equation}
where $S\subset [n]$, $S\ne\emptyset$;
$R$ ranges over all subsets of $[n]$ such that
$|R|=|S|$; $\sigma(R)$ and $\sigma(S)$ are the sums
of all elements of $R$ and $S$ respectively. Put
$$
S=\{
1,t+1,2t+1,\ldots,(t-1)t+1
\} .
$$
In the other words, to calculate $|B|$ we use the decomposition
by the first columns of the blocks.
Since  blocks of $B$ are regular, we see that if the
number of a row of the submatrix $B([n],S)$ does not
belong to $S$, then this row is zero.
Hence the sum in the right-hand side of (\ref{b9})
contains the term for $R=S$ only.
Because $|B(S,S)|=|b_{ij}^{(0)}|$, we have
$$
|B|=|b_{ij}^{(0)}|\cdot |B(\bar S,\bar S)|.
$$
The block partition of $B(\bar S,\bar S)$ is defined by
$\tilde\nu=(t-1,t-1,\ldots,t-1)$
and any block of $B(\bar S,\bar S)$ is regular.
Hence,
by the inductive assumption,
$$
|B(\bar S,\bar S)|=|b_{ij}^{(0)}|^{t-1}.
$$
Whence, (\ref{b9.0}) is proved.

To conclude the proof, it remains to note that,
in our case, $B_P=B_Q=(b_{ij}^{(0)})$
(see the proof of Theorem 2.3).

Suppose the partition $\nu$ is given by (\ref{b8.1}),
$B\in \bar {\cal B}_{\nu},$
$B_P=(B_{ij}^{\prime}),$
$B_Q=(B_{ij}^{\prime\prime}),$
i.e., $B_{ij}^{\prime}$ and $B_{ij}^{\prime\prime}$
are the blocks of $B_P$ and $B_Q$ respectively.

Recall that $B_P$ and $B_Q$ are block triangular and
respective diagonal blocks of these matrices
are coincide, i.e.,
$B_{ii}^{\prime}=B_{ii}^{\prime\prime},$
where $i=1,2,\ldots,l$. Then we have
\begin{equation}
|B_P|=|B_Q|=\prod _{i=1}^l
|B_{ii}^{\prime\prime}|.
\label{b10.0}
\end{equation}

{\bf Lemma 2.10.}
\begin{equation}
|B|=\prod _{i=1}^l
|B_{ii}^{\prime\prime}|^{\alpha_i}.
\label{b10}
\end{equation}
{\sc Proof.} Proceed by induction on the remainder
$(\alpha_1-\alpha_l)$.
If $\alpha_1-\alpha_l=0$, then (\ref{b10}) follows from
Lemma 2.9.

Suppose $(\alpha_1-\alpha_l)>0$. To calculate $|B|$ we apply
Laplace's theorem (see (\ref{b9})),
using the decomposition by the columns that contain
the first columns of the blocks $B_{i1},B_{i2},
\ldots,B_{is_1}$,
where $i\in [m]$, i.e., put
$$
S=\{ 1,\alpha_1+1,2\alpha_1+1,\ldots,(s_1-1)\alpha_1+1\}.
$$

Suppose $b_{pq}$ is an element of $B$ such that $q\in S$.
Since blocks of $B$ are regular and $\alpha_1>\alpha_2,$
we see that
if $p>\alpha_1s_1$, then $b_{pq}=0$. Moreover, $b_{pq}\ne 0$
implies $p\in S$.
Then with
$$
B(S,S)=B_{11}^{\prime}=B_{11}^{\prime\prime},
$$
from (\ref{b9}) it follows that
\begin{equation}
|B|=|B_{11}^{\prime\prime}|\cdot|\tilde B|,
\label{b11}
\end{equation}
where $\tilde B=B(\bar S,\bar S)=B_{S,S}.$

The block partition of $\tilde B$  is
defined by
$$
\tilde \nu=(
\underbrace{ \alpha_1-1,\ldots,\alpha_1-1}_{s_1},
\underbrace{\alpha_2,\ldots,\alpha_2}_{s_2},
\ldots,
\underbrace{\alpha_l,\ldots,\alpha_l}_{s_l}
).
$$
We claim that an arbitrary
block $\tilde B_{ij}$ of $\tilde B$ is regular.
Indeed, to obtain $\tilde B_{ij}$ from $B_{ij}$ we remove from
$B_{ij}$

i) the first row and the first column for $i,j\le s_1$;

ii) the first row for $i\le s_1$, $j>s_1$;

iii) the first column for $i>s_1$, $j\le s_1$.

Also, $\tilde B_{ij}=B_{ij}$ for $i>s_1$, $j>s_1$.
Since $B_{ij}$ is regular, we see that $\tilde B_{ij}$
is also regular.
Moreover, if $i\le s_1$, $j>s_1$, then
the diagonal  of the block $\tilde B_{ij}$ is zero.

Suppose the sets $\tilde P$, $\tilde Q$ are defined by
$J_{\tilde \nu}$ in the same way as $P,Q$ are defined by
$J_{\nu}$ (see above).
Denote by
$\tilde B_{ij}^{\prime\prime}$
the blocks of $\tilde B_Q$, i.e.,
$\tilde B_Q=(\tilde B_{ij}^{\prime\prime})$,
where $1\le i,j\le l+\sgn(\alpha_1-\alpha_2-1)-1$.

Suppose $\alpha_1-1>\alpha_2$.
Then, by the inductive assumption, we have
$$
|\tilde B|=|\tilde B_{11}^{\prime\prime}|^{\alpha_1-1}
\prod_{i=2}^l
|\tilde B_{ii}^{\prime\prime}|^{\alpha_i}.
$$
Since $B_{ii}^{\prime\prime}=\tilde B_{ii}^{\prime\prime}$,
where $i=1,2,\ldots,l$,
we see that (\ref{b10}) is proved.

Suppose $\alpha_1-1=\alpha_2$. Then
$$
\tilde \nu=(
\underbrace{ \alpha_2,\ldots,\alpha_2}_{s_1+s_2},
\underbrace{\alpha_3,\ldots,\alpha_3}_{s_3},
\ldots,
\underbrace{\alpha_l,\ldots,\alpha_l}_{s_l}
).
$$
By the inductive assumption, we have
\begin{equation}
|\tilde B|=
\prod_{i=1}^{l-1}
|\tilde B_{ii}^{\prime\prime}|^{\alpha_{i+1}}.
\label{b12}
\end{equation}
If $i\le s_1,$ $j>s_1$, then the diagonal of the block $\tilde B_{ij}$
of $\tilde B$ is zero. Hence $\tilde B_{11}^{\prime\prime}$
is block lower triangular and the diagonal blocks
of $\tilde B_{11}^{\prime\prime}$ are
$B_{11}^{\prime\prime}$ and $B_{22}^{\prime\prime}$.
Further,
$\tilde B_{ii}^{\prime\prime}=
B_{i+1,i+1}^{\prime\prime}$, where
$i=2,3,\ldots,l-1.$
Therefore from (\ref{b12}) it follows that
$$
|\tilde B|=
|B_{11}^{\prime\prime}|^{\alpha_2}
\cdot
|B_{22}^{\prime\prime}|^{\alpha_2}
\cdot
\prod_{i=3}^l
|B_{ii}^{\prime\prime}|^{\alpha_i}.
$$
Recall that $\alpha_2=\alpha_1-1$.
Then with (\ref{b11}) we arrive at (\ref{b10}).

{\sc Proof of  Proposition 2.2.} Follows immediately
from (\ref{b10.0}) and (\ref{b10}).

Let $m$ be a positive number. We say that $\mu=(\mu_1,\ldots,
\mu_l)$ is {\it a sequence} for $m$ if $\mu_1+\cdots+\mu_l=m$
and $\mu_1,\ldots,\mu_l$ are positive integers.
Suppose
the block partition of matrices of $M_m(K)$
is defined by the sequence $\mu$.

We say that
matrices  $C,\:\tilde C\in GL_m(K)$ are $\mu$-{\it mutual}
if

i) $C$ is lower block triangular;

ii) $\tilde C$ is upper block triangular;

iii)  respective diagonal blocks of $C$ and $\tilde C$ are equal.

If, in addition, all diagonal blocks of $C$ and $\tilde C$
are identity matrices, then $C$ and $\tilde C$ are
{\it unitary $\mu$-mutual}.

We say that matrices $A,\:B\in M_m(K)$ are $\mu$-{\it similar}
({\it unitary $\mu$-similar})
if there exist $\mu$-mutual
(respectively, unitary $\mu$-mutual)
matrices  \linebreak
$C,\tilde C\in GL_m(K)$
such that
$$
CA=B\tilde C.
$$
Evidently, if $\mu =(m)$, then the concept of the $\mu$-similarity
coincides with the concept of similarity in the classic case.

{\bf{Lemma 2.11.}} {\em{Matrices $A,\:B\in M_m(K)$ are $\mu$-similar
iff there exist an invertible lower block triangular matrix
$C=(C_{ij})$ and an upper block triangular matrix $D=(D_{ij})$ such that
$D_{ii}=C_{ii}^{-1}$ for $i=1,2,\ldots,l$ and the following
condition holds:
$$
CAD=B.
$$
}}

{\sc{Proof.}} Note that, for an arbitrary upper block triangular
matrix \linebreak
$G\in GL_m(K)$, the diagonal blocks of the matrix $H=G^{-1}$
are given by $H_{ii}=G_{ii}^{-1}$, where $i=1,2,\ldots,l.$  

Let $A\in M_n(K)$.

{\bf Theorem 2.4.} {\em If $B$ ranges over ${\cal B}_{\nu}$,
then $(BAB^{-1})_{P,Q}$ ranges over all matrices
that are $\mu$-similar to $A_{P,Q}$.}

{\sc Proof.} Follows from (\ref{b5}), Theorem 2.3 and Lemma 2.11.  

The following example shows the part of Theorem 2.4
in the solution of our main problem in the simplest case.

{\bf Example 2.1.} Consider a certain set of matrices that are similar
to a matrix of the form
$$
A=J_{\nu}+A_1\zeta,
$$
where $\nu=(2,1)$, $A_1\in M_3(K)$.

Let us prove that this set contains a unique matrix of the form
\begin{equation}
J_{\nu}+
\left(
\begin{array}{ccc}
\beta_{11}&0&0\\
\beta_{21}&0&\beta_{23}\\
\beta_{31}&0&\beta_{33}
\end{array}
\right)
\zeta ,
\label{e1}
\end{equation}
where
$
{\displaystyle
\left(
\begin{array}{cc}
\beta_{21}&\beta_{23}\\
\beta_{31}&\beta_{33}
\end{array}
\right)
}
\in {\cal R}_2,
$
${\cal R}_2$ is the set of canonical matrices
for  $(1,1)$-similar matrices
(see Example 3.1 and Theorem 3.1 below),
$\beta_{11}\in K$.

First let us prove that there exists a matrix $D$ of the form
(\ref{b2}) such that the matrix
$$
DAD^{-1}
$$
has  form (\ref{e1}). For all $B\in {\cal B}_{\nu}$,
we have
$$
BAB^{-1}=J_{\nu}+BA_1B^{-1}\zeta.
$$
By Theorem 2.4, $\{(BA_1B^{-1})_{P,Q} \; |\; B\in {\cal B}_{\nu}\}$
is the set of
all matrices \linebreak
that are $\mu$-similar
to $(A_1)_{P,Q}$, where $P=\{ 1\}$, $Q=\{ 2\}$, $\mu=(1,1)$. \linebreak
By Theorem 3.1, it follows that there exists a matrix
$B\in {\cal B}_{\nu}$ such that
$\tilde A_{P,Q} \in {\cal R}_2$,
where $\tilde A=BA_1B^{-1}$.

Further, consider the mapping
$$
J_{\nu}+\tilde A\zeta \mapsto
C(J_{\nu}+\tilde A\zeta)C^{-1},
$$
where $C=E_3+C_1\zeta$, $C_1\in M_3(K).$
The matrix $\tilde A_{P,Q}$ is unchanged by this mapping
(see Lemma 2.3). To be more precise,
$$
C(J_{\nu}+\tilde A\zeta)C^{-1}=
(E_3+C_1\zeta)(J_{\nu}+\tilde A\zeta)
(E_3-C_1\zeta)
$$
$$
=J_{\nu}+
\left(
\begin{array}{ccc}
\tilde a_{11}-c_{21}^{\prime}&
\tilde a_{12}+c_{11}^{\prime}-c_{22}^{\prime}&
\tilde a _{13}-c_{13}^{\prime}\\
\tilde a_{21}&
\tilde a_{22}+c_{21}^{\prime}&
\tilde a _{23}\\
\tilde a_{31}&
\tilde a_{32}+c_{31}^{\prime}&
\tilde a _{33}
\end{array}
\right)
\zeta,
$$
where $\tilde A=(\tilde a_{ij})$, $C_1=(c_{ij}^{\prime})$.
Now it is easy to see that there exists a matrix $C_1$
such that $C(J_{\nu}+\tilde A\zeta)C^{-1}$ has  form (\ref{e1}).

Finally, let us prove that if
$$
J_{\nu}+A_1^{\prime}\zeta \sim  J_{\nu}+A_1^{\prime\prime}\zeta
$$
and $A_1^{\prime},A_1^{\prime\prime}$ have  form (\ref{e1}),
then  $A_1^{\prime}=A_1^{\prime\prime}$.
Indeed, $(A_1^{\prime})_{P,Q}$,  $(A_1^{\prime\prime})_{P,Q}$
are $\mu$-similar and belong to ${\cal R}_2$.
Hence, by Theorem 3.1, we have
\begin{equation}
(A_1^{\prime})_{P,Q}=(A_1^{\prime\prime})_{P,Q}.
\label{e2}
\end{equation}
Since the trace of matrix is invariant with respect
to conjugation, from (\ref{e2}) it follows that
the element that belongs to the first row and to the first column
is the same for the matrices $ J_{\nu}+A_1^{\prime}\zeta$
and $J_{\nu}+A_1^{\prime\prime}\zeta$.

{\bf Remark.} It can be noted that the algorithm we use in the present section is similar in some aspects to this of Belitskii  (see, for example, \cite {Can}).

\section{ $\mu$-similar
matrices  for $\mu=(1,1,\ldots,1)$}
\setcounter{equation}{0}

In this section, suppose $\mu=(\underbrace{1,1,\ldots,1}_m)$,
where $m$ is a positive integer.

Let
$R,\:S$
be subsets of $[m]$ of the same cardinality. For
arbitrary $k,\:l\in [m]$,
by $R^k$and $S^l$ we denote the sets $R\cup \{ k\}$ and
$S\cup \{ l\}$ respectively.

Let $f_{R,S}$ be the mapping of $M_m(K)$ to $M_m(K)$ such that,
for a matrix $A\in M_m(K)$, an arbitrary element $y_{kl}$ of the matrix
$Y=f_{R,S}(A)$ is defined as follows:

i) if $k\notin R$ and $l\notin S$, then $y_{kl}=| A(R^k,S^l) |$;

ii) if $k\in R$ and $l\notin S$ or $k\notin R$ and $l\in S$,
then $y_{kl}=0$;

iii) $Y(R,S)=E$, where $E$ is the identity matrix.

Evidently, if $R=S=\emptyset$, then $f_{R,S}(A)=A$.
If $R=S=[m]$, then $f_{R,S}(A)=E_m$.

Let $A,C,D\in M_m(K)$;  $B=CAD$. Also put
$U= f_{R,R}(C)f_{R,S}(A)f_{S,S}(D)$,
$V= f_{R,S}(B)$, $X= f_{R,R}(C)$,
$Z= f_{S,S}(D)$ and as above  $Y=f_{R,S}(A)$.
Clearly, $u_{ij},\:v_{ij},\ldots,z_{ij}$ are the elements
of $U,\:V,\ldots,Z$ respectively, where $i,\:j\in [m]$.

{\bf {Lemma 3.1.}} {\em{ If $i\in R$ or $j\in S$,
then $u_{ij}=v_{ij}$.}}

{\sc Proof}. Suppose $i\in R$ and $j\notin S$. If $k\notin R$,
then $x_{ik}=0$. If $l\in S$, then $z_{lj}=0$. Therefore,
$$
u_{ij}=\sum_{\textstyle{{k\in R}\atop {l\notin S}}}
x_{ik}y_{kl}z_{lj}.
$$
But if $k\in R$ and $l\notin S$, then $y_{kl}=0$. Whence we have
$u_{ij}=0=v_{ij}$.

Similarly, if $i\notin R$ and $j\in S$, then
$u_{ij}=v_{ij}=0$.

Let $i\in R$ and $j\in S$. Since $Y(R,S)=E$, we have
$$
u_{ij}=\sum_{k\in R}x_{ik}z_{l(k)j},
$$
where $l(k)$ is the element of $S$ such that it has
the same ordinal number in the natural ordering of $S$
as the element $k\in R$ has. Recall that $(x_{ij})=
f_{R,R}(C)$ and $(z_{ij})=f_{S,S}(D)$. Then, by definition,
we have $x_{ik}=\delta_{ik}$ and $z_{l(k)j}=
\delta_{l(k)j}$, where $k\in R,\:l(k)\in S$. Thus we obtain
$$
u_{ij}=\sum_{k\in R} \delta_{ik}\delta_{l(k)j}=\delta_{ij}.
$$
Also, from the definition of $f_{R,S}$ it follows that
$v_{ij}=\delta_{ij}$.

Using the notations of Lemma 3.1, in addition
suppose:

i) the sets $R$ and $S$ are one-element, that is, $R=\{ r\}$
and $S=\{ s\}$, where $r,\:s\in [m]$;

ii) the matrix $C\;(D)$ is  lower (respectively, upper)
triangular and all off-diagonal elements of its $r$-th row
(respectively, $s$-th column) are equal to zero.

Then we have

{\bf{Lemma 3.2.}} {\em{If $i\notin R$ and $j\notin S$,
then $u_{ij}=v_{ij}$.}}

{\sc{Proof.}} This lemma is equivalent to the equality
\begin{equation}
\sum_{\textstyle{k\notin R \atop l\notin S}}
|C(R^i,R^k)|\cdot |A(R^k,S^l)| \cdot |D(S^l,S^j)|=
|B(R^i,S^j)|,
\label{c1}
\end{equation}
where $i\notin  R,\:j\notin S$. According to the Cauchy-Binet
formula we have
\begin{equation}
\sum_{T^{\prime},T^{\prime\prime}}
|C(R^i,T^{\prime})|\cdot |A(T^{\prime},T^{\prime\prime})|\cdot
|D(T^{\prime\prime},S^j)|=
|B(R^i,S^j)|,
\label{c2}
\end{equation}
where $T^{\prime},\:T^{\prime\prime}$ range independently
over all two-element subsets of $[m]$. Let us prove we can assume
that the sets $T^{\prime}$ and $T^{\prime\prime}$ of formula
(\ref{c2}) are only such that $r\in T^{\prime}$ and
$s\in T^{\prime\prime}$.
Indeed, let us check that if $r\notin T^{\prime}$
or $s\notin T^{\prime\prime}$, then
$$
|C(R^i,T^{\prime})|\cdot|A(T^{\prime},T^{\prime\prime})|
\cdot|D(T^{\prime\prime},S^j)|=0.
$$
By assumption, all off-diagonal elements of the $r$-th row
of the matrix $C$ are equal to zero. Therefore, if $r\notin T^{\prime}$,
then the matrix $C(R^i,T^{\prime})$ contains the zero row.
Hence $|C(R^i,T^{\prime})|=0$.

By analogy, the condition $s\notin T^{\prime\prime}$ yields
$|D(T^{\prime\prime},S^j)|=0$.  

It follows from Lemmas 3.1 and 3.2 that under the conditions
of   \linebreak
Lemma 3.2 we have the following result:

{\bf{Proposition 3.1.}}
$$
f_{R,S}(CAD)=f_{R,R}(C)f_{R,S}(A)f_{S,S}(D).  \; 
$$
%

An element $a_{rs}$ of a matrix $A=(a_{ij})$ is called
{\it marked} if $a_{rs}\ne 0$ and for all positive
integers $i,\:j$ such that $i\le r$, $j\le s$,
$(i,j)\ne (r,s)$, we have $a_{ij}=0$.
In other words, $a_{rs}$ is marked iff $a_{rs}$ is the only non-zero
element of the matrix $A([r],[s])$.
By ${\cal P}(A)$ denote the set of all marked element numbers of $A$.

For instance, if
$A=\left(
\begin{array}{ccc}
0&0&0\\
0&2&4\\
5&0&1
\end{array}
\right),
$
then the marked elements of $A$ are $a_{22}=2$ and $a_{31}=5$;
${\cal P}(A)=\{ (2,2),(3,1) \}$.
If
$A=\left(
\begin{array}{cc}
3&1\\
2&0\\
\end{array}
\right),
$
then the only marked element is $a_{11}=3$ and ${\cal P}(A)=\{
(1,1) \}$. Also, the matrix
$A=\left(
\begin{array}{cc}
0&0\\
0&0\\
\end{array}
\right)
$
has no marked elements, hence, ${\cal P}(A)=\emptyset$.

{\bf{Lemma 3.3.}} {\em{If matrices $A$ and $B$ are $\mu$-similar,
then ${\cal P}(A)={\cal P}(B)$.}}

{\sc{Proof.}} By Lemma 2.11, there exist  matrices \linebreak
$C=(c_{ij}),$ 
$D=(d_{ij})\in GL_m(K)$ such that

i) $C$ is lower triangular;

ii) $D$ is upper triangular;

iii) $d_{ii}=c_{ii}^{-1}$, where $i=1,2,\ldots,m$;

iv)  $B=CAD$.

 Suppose $(r,s)\in {\cal P}(A)$
and $i,\:j$ are positive integers such that $i\le r$,
$j\le s$ and $(i,j)\ne (r,s)$. We claim that $b_{ij}=0$
and $b_{rs}=c_{rr}c_{ss}^{-1}a_{rs}$.

In fact, since the matrices $C$ and $D$ are triangular,
we have
\begin{equation}
b_{ij}=\sum_{l=1}^j \sum_{k=1}^i c_{ik}a_{kl}d_{lj}.
\label{c4}
\end{equation}
The right-hand side of (\ref{c4}) contains just
the elements $a_{kl}$
of $A$ such that $k\le i$, $l\le j$,
and no others. But all these elements
are zero because $i\le r,$ $j\le s,$ $(i,j)\ne (r,s)$ and the element
$a_{rs}$ is marked. Therefore, $b_{ij}=0$.

In the same way, we get
\begin{equation}
b_{rs}=\sum_{k=1}^r \sum_{l=1}^s c_{rk}a_{kl}d_{ls}=
c_{rr}a_{rs}d_{ss}=c_{rr}a_{rs}c_{ss}^{-1},
\label{c5}
\end{equation}
as if $a_{kl}\ne 0$ for $k\le r$, $l\le s$, then $k=r$ and $l=s$.

Thus if $(r,s)\in {\cal P}(A),$ then $(r,s)\in {\cal P}(B)$,
that is, ${\cal P}(A)\subseteq {\cal P}(B)$.
Since the relation of  $\mu$-similarity is symmetric, we have
${\cal P}(B)\subseteq {\cal P}(A)$ and we see that the lemma
is proved.  

Note that the matrices $C$, $D^{-1}$ of the last proof
are $\mu$-mutual (see \linebreak
Section 2).

As above, suppose $a_{rs}$ is a marked element of $A\in M_m(K)$.

{\bf{Lemma 3.4.}} {\em{There exist $\mu$-mutual matrices $C$ and
$\tilde C$ such that the elements of the matrix $B=CAD$,
where $D=\tilde C^{-1}$,
satisfy the conditions

$b_{rj}=0$ for $j\ne s$ and $b_{is}=0$ for $i\ne r$
.}}

{\sc Proof.} Suppose $B$ is a matrix as required.
Then let us prove that the matrix equation $B=CAD$
is solvable with respect to the non-zero elements
of $C$ and $D$ as unknowns.

By definition, $c_{ij}=0$ for $i<j$ and $d_{kl}=0$
for $k>l$. Hence
for all $j\in \{ s+1,s+2,\ldots,m\}$, we have

\begin{equation}
b_{rj}=\sum_{k=1}^r \sum_{l=1}^j c_{rk}a_{kl}d_{lj}.
\label{c6}
\end{equation}
Since $a_{rs}$ is marked, we see that equality (\ref{c6})
can be written in the form
$$
b_{rj}=c_{rr}a_{rs}d_{sj}+ \sum_{k=1}^r \sum_{l=s+1}^j
c_{rk}a_{kl}d_{lj}.
$$
Because $a_{rs}\ne 0$
and $c_{rr}\ne 0$,
the condition
$b_{rj}=0$ determines uniquely the element $d_{sj}$, that is,
\begin{equation}
d_{sj}=-a_{rs}^{-1}c_{rr}^{-1}
\sum_{k=1}^r \sum_{l=s+1}^j c_{rk}a_{kl}d_{lj},
\label{c7}
\end{equation}
where $j=s+1,s+2,\ldots,m.$

By analogy, from $b_{is}=0$ it follows that

\begin{equation}
c_{ir}=-a_{rs}^{-1}d_{ss}^{-1}
\sum_{l=1}^s \sum_{k=r+1}^i c_{ik}a_{kl}d_{ls},
\label{c8}
\end{equation}
where $i=r+1,r+2,\ldots,m$.

By $Q_1$ (respectively, $Q_2$) denote the set of the elements
$c_{ij},\:d_{kl}$ that enter in the notations of the right
(respectively, left)-hand side of  equalities (\ref{c7})
and (\ref{c8}), that is,
$$
Q_1=\{ c_{rk}|k\le r \}\cup \{c_{ik}|r<k\le i\}
\cup \{d_{ls}|l\le s\}\cup \{ d_{lj}|s<l\le j\},
$$
$$
Q_2=\{ c_{ir}|i>r\} \cup \{ d_{sj}|j>s\}.
$$

We have $Q_1\cap Q_2=\emptyset$.
Further, from Lemma 3.3 it follows that
$b_{rj}=0$ for $j<s$ and $b_{is}=0$ for $i<r$.
Thus we see that the lemma is
proved.

In the following suppose the marked element $a_{rs}$ is minimal in
${\cal P}(A)$ with respect to the ordering by the numbers
of rows of  $A$.

By definition, put
$$
Q_3=\{ c_{rk}|k<r\} \cup \{ d_{ls}|l<s\}.
$$
Evidently, $Q_3\subseteq Q_1$.

An element of $Q_i$ is called a {\it $Q_i$-element},
 where $i=1,2,3.$

Let $A,B,C,D$ be the matrices of Lemma 3.4.
Every element of the matrix $B=CAD$ is a polynomial in 
elements of $C,A,D$. For all these polynomials, using
equalities
(\ref{c7}) and (\ref{c8}), express all $Q_2$-elements
in terms of the  $Q_1$-elements. Then we have

{\bf Lemma 3.5.} {\em Every element of the matrix $B$ is independent of
$Q_3$-elements.}

{\sc Proof.} By Lemma 3.3, if $i\le r$, $j\le s$ and $(i,j)\ne (r,s)$, then
$b_{ij}=0$. Also,
by assumption, we have $b_{rj}=0$ for $j>s$ and $b_{is}=0$
for $i>r$. Further, $b_{rs}=c_{rr}c_{ss}^{-1}a_{rs}$
(see (\ref{c5})). Whence suppose
$i>r$ or $j>s$, and $(i,j)\ne (r,s)$.

We have

\begin{equation}
b_{ij}=\sum_{k=1}^i\sum_{l=1}^j c_{ik}a_{kl}d_{lj}.
\label{c9}
\end{equation}
This notation alone does not contain any $Q_3$-elements.
Hence  $Q_3$-elements do not appear in the notation of
$b_{ij}$ until we write the
$Q_2$-elements in terms of  $Q_1$-elements.

Let $i>r$ and $j<s$. Then
 the right-hand side of (\ref{c9})
does not contain any $Q_2$-elements of  $D$.
The monomials of (\ref{c9}) that contain $Q_2$-elements
of $C$ are of the form $c_{ir}a_{rl}d_{lj}$, where $1\le l\le j$.
But all these monomials are equal to zero, since $a_{rl}=0$
for $l\le j<s$.

Similarly, the case when $i<r$ and $j>s$ can be considered.

Let $i>r$ and $j>s$. Write the element $b_{ij}$ in the form
$b_{ij}=b_{ij}^{\prime}+b_{ij}^{\prime\prime}$,
 where the notation of $b^{\prime}_{ij}$ does not
contain any $Q_2$-element and every monomial of
$b^{\prime\prime}_{ij}$
contains some $Q_2$-element. To be precise,
taking into account that $a_{rs}$ is marked, we have
$$
b^{\prime}_{ij}=\sum_{\textstyle {k=1 \atop k\ne r}}^i
 \sum_{\textstyle{l=1 \atop l\ne s}}^j
c_{ik}a_{kl}d_{lj},
$$
\begin{equation}
b^{\prime\prime}_{ij}=d_{sj}\sum_{k=r+1}^ic_{ik}a_{ks}+
c_{ir}\sum_{l=s+1}^j a_{rl}d_{lj}+c_{ir}a_{rs}d_{sj},
\label{c10}
\end{equation}
where $c_{ir},\:d_{sj}\in Q_2$ and the other elements of the
right-hand side of (\ref{c10}) do not
belong to $Q_2$. Change $c_{ir}$ and $d_{sj}$ by  equalities
(\ref{c7}) and (\ref{c8}) respectively.
Then we obtain
$$
b^{\prime\prime}_{ij}=-a_{rs}^{-1}c_{rr}^{-1}
\sum_{t=1}^r \sum_{l=s+1}^j \sum_{k=r+1}^i c_{rt}a_{tl}
d_{lj}c_{ik}a_{ks}
$$
$$
-a_{rs}^{-1}d_{ss}^{-1}\sum_{q=1}^s \sum_{k=r+1}^i
\sum_{l=s+1}^j c_{ik}a_{kq}d_{qs}a_{rl}d_{lj}
$$
\begin{equation}
+a_{rs}^{-1} c_{rr}^{-1} d_{ss}^{-1}
\sum_{t=1}^r \sum_{l=s+1}^j \sum_{q=1}^s \sum_{k=r+1}^i
c_{rt}a_{tl}d_{lj}c_{ik}a_{kq}d_{qs}.
\label{c11}
\end{equation}
By $\sigma_{rt}$ denote the coefficient of $c_{rt}$ in the right-hand
side of (\ref{c11}), where $t<r$. Then we have
$$
\sigma_{rt}=-a_{rs}^{-1}c_{rr}^{-1}
\sum_{l=s+1}^j \sum_{k=r+1}^i a_{tl}d_{lj}c_{ik}a_{ks}
$$
$$
+a_{rs}^{-1}c_{rr}^{-1}d_{ss}^{-1}
\sum_{l=s+1}^j \sum_{q=1}^s \sum_{k=r+1}^i
a_{tl}d_{lj}c_{ik}a_{kq}d_{qs}.
$$
Since the marked element $a_{rs}$ is minimal, we have
$a_{tl}=0$ for $t<r$. Hence $\sigma_{rt}=0$, where
$t=1,2,\ldots,r-1.$

By $\tau_{qs}$ denote the coefficient of $d_{qs}$ in
the right-hand side of (\ref{c11}), where $q<s$. Then
we see that
$$
\tau_{qs}=-a_{rs}^{-1}d_{ss}^{-1}
\sum_{k=r+1}^i \sum_{l=s+1}^j c_{ik}a_{kq}a_{rl}d_{lj}+
a_{rs}^{-1}c_{rr}^{-1}d_{ss}^{-1}
\sum_{t=1}^r \sum_{l=s+1}^j \sum_{k=r+1}^i
c_{rt}a_{tl}d_{lj}c_{ik}a_{kq}.
$$

With $a_{tl}=0$ for $t<r$ we get
$$
\tau_{qs}=-a_{rs}^{-1}d_{ss}^{-1}
\sum_{k=r+1}^i \sum_{l=s+1}^j c_{ik}a_{kq}a_{rl}d_{lj}+
a_{rs}^{-1}d_{ss}^{-1}
\sum_{l=s+1}^j \sum_{k=r+1}^i
a_{rl}d_{lj}c_{ik}a_{kq}\equiv 0,
$$
where $q=1,2,\ldots,s-1$. We see that the lemma is proved. 

Note that from this proof it follows that Lemma 3.5 also holds
if the element $a_{rs}\in {\cal P}(A)$ is minimal with respect to
the ordering by the numbers of columns.

As above, suppose $A\in M_m(K),\:a_{rs}\in {\cal P}(A)$
and $a_{rs}$ is minimal with respect to the ordering
by the numbers of rows. Let $C$ and $D$ be the matrices
of Lemma 3.4.

{\bf Lemma 3.6.} {\em It can be assumed that all off-diagonal
elements of the $r$-th row of $C$ and of the $s$-th column
of $D$ are equal to zero.}

{\sc Proof.} Let $Q_3$ be as above.
Recall that $Q_3\subseteq Q_1$.
By Lemma 3.5, every element of $B$
does not depend on $Q_3$-elements. Hence all elements
of $Q_3$ can be equated to zero without violating the equality
$B=CAD$ and the lemma is proved.

Let us remark that Lemmas 3.4-3.6  are correct as before
if, in addition, we assume that the matrices $C$ and $D^{-1}$
are  unitary $\mu$-mutual. Indeed, suffice it to note
that  diagonal elements of $C$ and $D$ do not belong to the
set $Q_2\cup Q_3$ of `dependent parameters'. Hence we can assume
$c_{ii}=1$, where $i=1,2,\ldots,m.$


By definition, put
$$
H_{t+1}=\diag(\underbrace{-1,\ldots,-1}_t,
\underbrace{1,\ldots,1}_{m-t-1}),
$$
where $t=0,1,2,\ldots,m-1$. Evidently, $H_1=E_{m-1}$ and
$H_m=-E_{m-1}$.

Under the conditions of Lemma 3.6, in addition, suppose
all diagonal elements of the matrices $C$ and $D$ are equal to 1
(in other words, $C$ and $D^{-1}$ are unitary $\mu$-mutual).
Also, by definition, put
$$
A_{\langle R,S \rangle}=
Y_{R,S}H_s a_{rs}^{-1},
$$
where $Y=f_{R,S}(A)$, $R=\{ r\}$, $S=\{ s \}$ . Then we have

{\bf Proposition 3.2.}

\begin{equation}
B_{R,S}=C_{R}A_{\langle R,S\rangle}D_{S}.
\label{c12}
\end{equation}

{\sc Proof.} From Proposition 3.1 and Lemma 3.6 it follows that
$$
f_{R,S}(B)=f_{R,R}(C)Yf_{S,S}(D).
$$
By definition, excepting the diagonal element,
the $r$-th column of $f_{R,R}(C)$ is  zero.
Hence  condition $(R)$ holds
for the matrix $f_{R,R}(C)$ (see Section 2).

By analogy,  condition $(\bar S)$ holds for the matrix
$f_{S,S}(D)$.
Therefore, from (\ref{b4}) it follows that
\begin{equation}
(f_{R,S}(B))_{R,S}=(f_{R,R}(C))_{R}Y_{R,S}
(f_{S,S}(D))_{S}.
\label{c13}
\end{equation}

Put $X=f_{R,R}(C)$. Then, under the accepted conditions for $C$,
we see that $x_{ij}=-c_{ij}$, where $i>r$ and $j<r$;
$x_{ir}=0,$ where $i\ne r$; otherwise $x_{ij}=c_{ij}$.
Whence,
\begin{equation}
(f_{R,R}(C))_{R}=H_rC_{R}H_r.
\label{c14}
\end{equation}
Likewise,
\begin{equation}
(f_{S,S}(D))_{S}=H_sD_{S}H_s.
\label{c15}
\end{equation}
Using (\ref{c13})-(\ref{c15}), we get
$$
H_r(f_{R,S}(B))_{R,S}H_s=C_{R}H_rY_{R,S}H_sD_{S}.
$$
Recall that $b_{rj}=0$ for $j\ne s$ and $b_{is}=0$ for $i\ne r$;
$b_{rs}=a_{rs}$ and $b_{ij}=0$ for $i<r$.

Put $V=f_{R,S}(B)$. We see that
$v_{ij}=0$ for $i<r$;
$v_{rj}=0$ for $j\ne s$.
Also,
$v_{is}=0$ for $i\ne s$;
$v_{rs}=1$.
Further,
if $i>r$ and $j<s$, then $v_{ij}=-b_{ij}a_{rs}$;
if $i>r$ and $j>s$, then $v_{ij}=b_{ij}a_{rs}$.
This yields that

$$
H_r(f_{R,S}(B))_{R,S}H_s=a_{rs}B_{R,S}.
$$
Inasmuch as the element $a_{rs}$ is minimal, the first
$(r-1)$ rows of $A$ are zero. Hence the first $(r-1)$ rows
of $Y=f_{R,S}(A)$ are also zero. Therefore,
$$
H_rY_{R,S}=Y_{R,S}. \; 
$$
We see that the proposition is proved.

By ${\cal U}_m$ denote the subset of $M_m(K)$ such that
$A\in M_m(K)$ belongs to ${\cal U}_m$ iff any row and any column of $A$
contain not more than one non-zero element.

{\bf Lemma 3.7.} {\em For any $A\in M_m(K)$ there exists a matrix
$B\in {\cal U}_m$ such that the matrices $A$ and $B$
are unitary $\mu$-similar.}

{\sc Proof.} The proof is by induction on $m$. For $m=1$,
there is nothing to prove.

Suppose $m>1$, $A\in M_m(K)$ and $a_{rs}$ is the minimal marked element
of $A$.

In (\ref{c12}), the matrix $A_{\langle R,S\rangle}$ is of
order $(m-1)$. $C_R$ and $D_S$ are unitary $\mu$-mutual.
If an element of the matrix $C_R$ ($D_S$) is
off-diagonal and non-zero then it does
not belong to the set $Q_2\cup Q_3$.
Therefore, $C_R-E_{m-1}$ ($D_S-E_{m-1}$) is the strictly
lower (respectively, upper) triangular matrix of the general form.
With the inductive assumption,
the lemma
is proved.

{\bf Lemma 3.8.} {\em Suppose matrices $A,B\in {\cal U}_m$
are unitary $\mu$-similar. Then
$$
A=B.
$$
}
{\sc Proof.} Let us prove that if $A\in {\cal U}_m$,
then
\begin{equation}
A_{\langle R,S\rangle}=A_{R,S},
\label{c17}
\end{equation}
where $R=\{ r\},$ $S=\{ s\}$ and $(r,s)$ is the number of the
minimal marked element of $A$.

Indeed, let $Y=f_{R,S}(A)$. By definition, $y_{rs}=1$;
$y_{rj}=0$, where $j\ne s$; $y_{is}=0$, where $i\ne r$.

Suppose $i,j\in [m]$ are integers such that $i\ne r$,
$j\ne s$ and $a_{ij}=0$. Since $A\in {\cal U}_m$, we have
$a_{is}=0$. Therefore, $y_{ij}=0$.

Suppose $i,j\in [m]$ are integers such that $i\ne r$,
$j\ne s$ and $a_{ij}\ne 0$.
As above, $A\in {\cal U}_m$ implies $a_{is}=0$.
Hence, $y_{ij}=a_{ij}a_{rs}\sgn (j-s).$

Now with the definition of $A_{\langle R,S\rangle}$,
we arrive at (\ref{c17}).

The proof of the lemma is by induction on $m$.
The case $m=1$ is trivial. Suppose $m>1$ and $B=CAD$,
where $C,D^{-1}$ are unitary $\mu$-mutual.
By Lemma 3.3, with (\ref{c5}) we see that
the minimal marked element of $B$
is $b_{rs}=a_{rs}$. From (\ref{c12}) and (\ref{c17})
it follows that
$$
B_{R,S}=C_RA_{R,S} D_S.
$$
Since $A,B\in {\cal U}_m$, we have $A_{R,S},B_{R,S}\in {\cal U}_m$.
Evidently, $C_R$ and $D_S^{-1}$ are unitary $\mu$-mutual.
Whence, by the inductive assumption, $A_{R,S}=B_{R,S}$.
Therefore, $A=B$ and the lemma is proved.

{\bf Proposition 3.3.} {\em ${\cal U}_m$ is a set of
canonical matrices
for unitary $\mu$-similar
matrices of order $m$,
where $\mu=(\underbrace{1,\ldots,1}_m)$.}

{\sc Proof.} Follows immediately from Lemmas 3.7 and 3.8.

Let $A\in {\cal U}_m$ and
\begin{equation}
\langle
a_{i_1j_1},a_{i_2j_2},\ldots,a_{i_lj_l}
\rangle
\label{x.1}
\end{equation}
be a sequence of non-zero elements of $A$.
We say that this sequence is {\it a chain} of $A$ if
$l=1$ or $i_1<i_{t+1}$ and $i_{t+1}=j_t$
for $t=1,2,\ldots,l-1.$ In this case, $a_{i_1j_1}$
is {\it the leading element} of the chain.
A chain is {\it maximal} if it can not be imbeded in another chain
of larger lenght. In the following we consider only
maximal chains. A chain of lenght $l$ is
$closed$ if $i_1=j_l$ (see (\ref{x.1})); otherwise it is {\it open}.
For example, if
$
A=
\left(
\begin{array}{ccc}
0&0&0\\
0&0&1\\
2&0&0
\end{array}
\right),
$
then $\langle a_{23},a_{31} \rangle$ is the open chain.
If 
$
A=
\left(
\begin{array}{ccc}
0&5&0\\
0&0&3\\
4&0&0
\end{array}
\right),
$
then $\langle a_{12},a_{23},a_{31} \rangle$ is the closed chain.
If
$
A=
\left(
\begin{array}{ccc}
0&6&0\\
7&0&0\\
0&0&8
\end{array}
\right),
$
then $A$ contains two closed chains:
$\langle a_{12}, a_{21} \rangle$ and
$\langle a_{33} \rangle$.
It follows from the definition of ${\cal U}_m$,
that if two chains of a matrix have a common element
then these chains coincide.

By ${\cal R}_m$ denote the subset of ${\cal U}_m$
such that $A\in {\cal U}_m$
belongs to ${\cal R}_m$ iff,
possibly except leading elements of closed chains,
all non-zero elements of $A$ are units.

{\bf Example 3.1.} Let $\mu =(1,1)$.
Then we have
$$
{\cal R}_2 =
\left\{
\left(
\begin{array}{cc}
\alpha&0\\
0&\beta
\end{array}
\right)
\;
\alpha,\beta\in K;
\;
\left(
\begin{array}{cc}
0&\gamma\\
1&0
\end{array}
\right)
\;
\gamma\in K;
\;
\left(
\begin{array}{cc}
0&1\\
0&0
\end{array}
\right)
\right\}. \; 
$$

{\bf Theorem 3.1.} {\em ${\cal R}_m$ is a set of canonical
matrices
for  $\mu$-similar matrices
of order $m$, where $\mu=(\underbrace{1,\ldots,1}_m)$.}

To prove this theorem, we need several lemmas.

{\bf Lemma 3.9.} {\em For any $A\in M_m(K)$ there exists a matrix
$B\in {\cal R}_m$ such that $A,\:B$ are $\mu$-similar.}

{\sc Proof.} By Proposition 3.3, there exists a unique matrix
$A^{\prime}\in {\cal U}_m$ such that $A$ and $A^{\prime}$
are unitary $\mu$-similar. Therefore, it suffices to show
that there exists a diagonal matrix
$C=\diag(\gamma_1,\gamma_2,\ldots,\gamma_m)\in GL_m(K)$
such that $B=CA^{\prime}C^{-1}\in {\cal R}_m.$

In the mapping $A'\mapsto B=CA'C^{-1}$, elements of different chains
are transformed independently.
Let $\langle a'_{i_1i_2}, a'_{i_2i_3}, a'_{i_3i_4},
\ldots,a'_{i_li_{l+1}} \rangle$ be a chain of $A'$.
Put
$\gamma_{i_{l+1}}=1,$
\begin{equation}
\gamma_{i_t}=\prod_{k=t}^{l}(a'_{i_k i_{k+1}})^{-1},
\label{x.2}
\end{equation}
where $t=2,3,\ldots,l$.
Then we have $b_{i_ti_{t+1}}=
\gamma_{i_t}\gamma_{i_{t+1}}^{-1}a'_{i_t i_{t+1}}=1$,
where $t=2,3,\ldots,l$.

If the chain is open, i.e. $i_{l+1}\ne i_1$,
define $\gamma_{i_1}$ by (\ref{x.2}) for $t=1$.
Then we obtain
$b_{i_1i_2}=\gamma_{i_1}\gamma_{i_2}^{-1} a'_{i_1i_2}=1$.

If the chain is closed, i.e. $i_{l+1}=i_1$, we have
$b_{i_1i_2}=\gamma_{i_1}\gamma_{i_2}^{-1} a'_{i_1i_2}=
\gamma_{i_2}^{-1} a'_{i_1 i_2}=
a'_{i_1i_2} a'_{i_2i_3}\cdots a'_{i_li_1}.$ 

Thus we obtain the diagonal matrix $C$, such that
$CA^{\prime}C^{-1}\in {\cal R}_m$ and the lemma is proved.

{\bf Lemma 3.10.} {\em If matrices $A,B\in {\cal R}_m$
are $\mu$-similar, then $A=B$.}

{\sc Proof.} By Lemma 2.11, there exist
matrices $C,D\in GL_m(K)$
such that $C$ is lower triangular, $D$ is upper triangular and
$$
B=CAD.
$$
The matrix $C$ (respectively, $D$)
can be uniquely represented in the form
$C=C_1C_2$ ($D=D_2D_1$), where $C_2-E_m$ ($D_2-E_m$)
is strictly lower (upper) triangular,
$C_1$ $(D_1)$ is invertible and diagonal.
By Lemma 2.11, $D_1=C_1^{-1}$. Hence we have
$$
B=C_1C_2AD_2C_1^{-1}.
$$
Since ${\cal R}_m\subseteq {\cal U}_m$, we have
$A,B\in {\cal U}_m$. Any conjugation by a diagonal matrix
preserves the set of  numbers of non-zero elements
of the matrix.
Hence, $C_1^{-1}BC_1\in {\cal U}_m$.
This yields that
$C_1^{-1}BC_1$ and $A$ are unitary $\mu$-similar.
Whence, by Lemma 3.8, we have $A=C_1^{-1}BC_1$ or
\begin{equation}
B=C_1AC_1^{-1}.
\label{c18}
\end{equation}

By the definition of ${\cal R}_m$, possibly except leading
elements of closed chains,
all non-zero elements of $A$ and $B$ are units.

Suppose
$\langle a_{i_1i_2},a_{i_2i_3},\ldots,a_{i_li_1}
\rangle$,
$\langle b_{i_1i_2},b_{i_2i_3},\ldots,b_{i_li_1}
\rangle$, are respective closed chains of $A$ and $B$.
Since $a_{i_1i_2}\cdots a_{i_li_1}=
b_{i_1i_2}\cdots b_{i_li_1}$ and
$a_{i_ti_{t+1}}=b_{i_ti_{t+1}}$ for $t=2,3,\ldots,l,$
where $i_{l+1}=i_1$, we have $a_{i_1i_2}=b_{i_1i_2}$. Thus the lemma is proved.

{\sc Proof of Theorem 3.1.} Follows immediately from Lemmas 3.9
and 3.10.

\section{Some special classes of similar matrices }
\setcounter{equation}{0}

In the present section we find a canonical matrix
for every class of matrices
that are similar to an arbitrary matrix of the form
$J_{n,\alpha}+A_1\zeta$, where $\alpha\in K$,
$A_1\in M_n(K)$.

Let
\begin{equation}
A=J_{n,0}+A_1\zeta,
\label{d1}
\end{equation}
where $A_1\in M_n(K)$. By $p_i$ and $q_j$ denote the traces
 the matrices $A^i$ and $A_1J_{n,0}^{j-1}\zeta$
respectively, i.e.,
$$
p_i=\tr A^i,
$$
\begin{equation}
q_j=\tr A_1J_{n,0}^{j-1}\zeta,
\label{d2}
\end{equation}
where $i=1,2,\ldots;\:j=2,3,\ldots. $
Also, by definition, put
$p_1=q_1.$

By $a_{ij}^{\prime}$ denote   elements of $A_1$.

{\bf Lemma 4.1.} {\em
$$
q_k=\sum_{i-j=k-1}a_{ij}=\sum_{i-j=k-1}a_{ij}^{\prime}\zeta,
$$
where $k=1,2,\ldots,n.$
}

{\sc Proof.} As above,  matrix units are denoted by $e_{ij}$.
Since
${\displaystyle
J_{n,0}=\sum_{i=1}^{n-1}e_{i,i+1},
}$
we have \\
${\displaystyle
J_{n,0}^k=\sum_{i=1}^{n-k}e_{i,k+i}.
}$
Because
${\displaystyle
A_1=\sum_{i=1}^n \sum_{j=1}^n a_{ij}^{\prime}e_{ij},
}$
we obtain
$$
A_1J_{n,0}^k=\sum_{i=1}^n \sum_{j=1}^n
a_{ij}^{\prime}e_{ij}
\sum_{t=1}^{n-k}e_{t,k+t}=
\sum_{i=1}^n \sum_{j=1}^{n-k} a_{ij}^{\prime}e_{i,k+j}.
$$
Therefore,
$$
\tr A_1 J_{n,0}^k=\sum_{i-j=k}a_{ij}^{\prime}.
$$
With (\ref{d2})
we obtain
${\displaystyle
q_k=\sum_{i-j=k-1}
a_{ij}^{\prime}\zeta}$,
where $k\in [n]$.
If an element $x_{ij}$ of $J_{n,0}$
is not zero, then $i<j$.
Hence, for an arbitrary $k\in [n]$, we have
$$
\sum_{i-j=k-1}
a_{ij}=
\sum_{i-j=k-1}
a_{ij}^{\prime}\zeta. \;  
$$
Thus the lemma is proved.

For an arbitrary positive integer $k$, we have
$$
A^k=(J_{n,0}+A_1\zeta)^k=J_{n,0}^k+
\sum_{i=0}^{k-1}J_{n,0}^iA_1J_{n,0}^{k-i-1}\zeta.
$$
Hence,
$$
p_k=\tr A^k=\tr J_{n,0}^k+
\sum_{i=0}^{k-1} \tr J_{n,0}^iA_1 J_{n,0}^{k-i-1}\zeta=
k\tr A_1 J_{n,0}^{k-1}\zeta.
$$
Now from definition (\ref{d2}) it follows that
\begin{equation}
p_k=kq_k,
\label{d3}
\end{equation}
where $k=1,2,\ldots$.

Let $A$ be given by  (\ref{d1}),
$B=J_{n,0}+B_1\zeta$,
where $B_1\in M_n(K)$.

{\bf Lemma 4.2.} {\em If the matrices $A,\:B$ are similar, then
$$
q_k(A)=q_k(B),
$$
where $k=1,2,\ldots.$ }

{\sc Proof.} Since the function $\tr$ is invariant with
respect to conjugation, we see that this lemma
follows from (\ref{d3}).  

As above, $A=J_{n,0}+A_1\zeta,$ where $A_1=(a_{ij}^{\prime})
\in M_n(K)$.

{\bf Lemma 4.3.} {\em There exists a unique matrix $B_1
\in M_n(K)$ such that

i) $A$ is similar to $B=J_{n,0}+
B_1\zeta;$

ii) all columns of $B_1$, from the second one,
are zero.}

{\sc Proof.}  From Lemmas 4.1 and 4.2 it follows that if  conditions
i) and ii) hold for some matrix $B_1\in M_n(K)$, then
$B_1$ is unique. Therefore, let us prove that a required
matrix $B_1$   exists.

Let $C=E_n+C_1\zeta,$ where $C_1=(c_{ij}^{\prime})\in M_n(K)$.
Then $C^{-1}=E_n-C_1\zeta$ and we have
$$
CAC^{-1}=(E_n+C_1\zeta)(J_{n,0}+A_1\zeta)(E_n-C_1\zeta)=
J_{n,0}+(A_1+C_1J_{n,0}-J_{n,0}C_1)\zeta.
$$
By $b_{ij}^{\prime}$ denote  elements of the matrix
$B_1=A_1+C_1J_{n,0}-J_{n,0}C_1$.
We see that
\begin{equation}
b_{ij}^{\prime}=
a_{ij}^{\prime}+d_{i,j-1}-d_{i+1,j},
\label{d4}
\end{equation}
where
$$
d_{uv}=\left\{
\begin{array}{l}
c_{uv}^{\prime} \qquad
\mbox{if}\;\;
(u,v)\in [n]\times [n],   \\
0 \qquad
\mbox{otherwise.}
\end{array}
\right.
$$
Suppose $i,j,k,l\in [n]$ are integers such that $i-j\ne k-l$.
Then $b_{ij}^{\prime}$ and $b_{kl}^{\prime}$
do not contain the same elements of $C_1$.

Consider the set of the equalities $b_{ij}^{\prime}=0$,
where $i-j=\mbox{const} \ge 0$ and $j>1$, as a system of  equations
with respect to $c_{ts}^{\prime}$ as unknowns.
The matrix of this system is triangular with  units
by the diagonal.
Hence this system has a unique
solution.

Now suppose  $i-j=\mbox{const}<  0$. Then $b_{ij}^{\prime}=0$
is the system of $i-j+n$ linear equations with $i-j+n+1$
unknowns. The matrix of this system contains the triangular
submatrix of order $i-j+n$ with  units by the diagonal.
Hence this system is solvable and the lemma is proved.

The last lemma yields the following result.

Suppose $A=A_0+A_1\zeta$, where $A_0,A_1\in M_n(K)$,
$A_0 \stackrel{K}{\sim} J_{n,\alpha}$,  $\alpha\in K$.
Then we have

{\bf Theorem 4.1.} {\em The class of all matrices that
are similar to $A$ contains a unique matrix of the form
$$
J_{n,\alpha}+B_1\zeta,
$$
where
$$
B_1=
\left(
\begin{array}{cccc}
\beta_1&0&\cdots&0\\
\beta_2&0&\cdots&0\\
\cdots&\cdots&\cdots&\cdots\\
\beta_n&0&\cdots&0
\end{array}
\right)
,
$$
$\beta_1,\beta_2,\ldots,\beta_n\in K$.  
}

Using the results of the present paper,
we produce canonical matrices for $M_n({\cal D})$,
where $n=2,3.$

First consider the following easy obtainable result.
Suppose $A=\alpha E_n+A_1\zeta$,
where $\alpha \in K$, $A_1\in M_n(K)$;
$J\in M_n(K)$
is a matrix such that it
 has a Jordan normal form and $J\stackrel{K}{\sim} A_1$.

{\bf Lemma 4.4.} {\em The class of the matrices that are
similar to $A$ contains a  matrix of the form
$$
\alpha E_n+J\zeta
$$
and the matrix $J$ is determined uniquely up to permutation
of Jordan blocks of equal orders.

}

{\sc Proof.} 
Suppose
$C=C_0+C_1\zeta$, where 
$C_0,C_1\in M_n(K)$ and $C$ 
is invertible.
Then we have
$C^{-1}=C_0^{-1}(E_n-C_1C_0^{-1}\zeta).$ Hence,
$$
C(\alpha E_n+B_1\zeta)C^{-1}=
\alpha E_n+C_0B_1C_0^{-1}\zeta.
$$
Now suffice it to apply the Jordan theorem.

{\bf Example 4.1.} The following matrices are canonical for
$M_2({\cal D})$:

i) $\diag(\alpha_1+\beta_1\zeta,\alpha_2+\beta_2\zeta)$
(up to permutation of  diagonal elements),
where $\alpha_i,\beta_j\in K,$ $\alpha_1\ne\alpha_2$
(see Theorem 2.2);

ii) $J_{2,\alpha_3}+
\left(
\begin{array}{ll}
\beta_3&0\\
\beta_4&0
\end{array}
\right)
\zeta,
$
where $\alpha_3, \beta_3,\beta_4 \in K$ (see Theorem 4.1);

iii) $\alpha_4E_2+J\zeta$, where $\alpha_4\in K$ and
$J\in M_2(K)$ is of the form: $\diag(\gamma_1,\gamma_2)$
(up to permutation of diagonal elements) or
$J_{2,\gamma_3}$, where $\gamma_i\in K$
(see \linebreak
Lemma 4.4).

{\bf Example 4.2.} The following matrices are canonical for
$M_3({\cal D})$:

i) $\diag(\alpha_1+\beta_1\zeta,\alpha_2+\beta_2\zeta,
\alpha_3+\beta_3\zeta)$ (up to permutation of  diagonal  elements),
where $\alpha_i,\beta_j\in K$, $\alpha_1,\alpha_2,\alpha_3$
are pairwise distinct (see  \linebreak
Theorem 2.2);

ii) $A\dot + E_1(\alpha_4+\beta_4\zeta)$,
where $A$ is a canonical matrix of $M_2({\cal D})$
(see \linebreak
Example 4.1), $\alpha_4,\beta_4\in K$ and $\alpha_4$
is not a root of the characteristic polynomial of the matrix
$A|_{\zeta=0}\in M_2(K)$ (see Theorem 2.2);

iii) $(J_{2,\alpha_5}\dot +J_{1,\alpha_5})+
\left(
\begin{array}{lll}
\beta_{11}&0&0\\
\beta_{21}&0&\beta_{23}\\
\beta_{31}&0&\beta_{33}
\end{array}
\right)
\zeta,$
where $\alpha_5,\beta_{ij}\in K$,\\
$
\left(
\begin{array}{ll}
\beta_{21}&\beta_{23}\\
\beta_{31}&\beta_{33}
\end{array}
\right)
\in {\cal R}_2
$
(see Examples 2.1 and 3.1);

iv) $J_{3,\alpha_6}+
\left(
\begin{array}{lll}
\beta_5&0&0\\
\beta_6&0&0\\
\beta_7&0&0
\end{array}
\right)
\zeta,$
where $\alpha_6,\beta_i\in K$ (see Theorem 4.1);

v) $\alpha_7E_3+J\zeta$, where $J\in M_3(K)$ is of the form:
$\diag(\gamma_1,\gamma_2,\gamma_3)$ (up to permutation of  diagonal
elements) or $J_{2,\gamma_4}\dot +J_{1,\gamma_5}$ or $J_{3,\gamma_6}$,
$\gamma_i\in K$ (see \linebreak
Lemma 4.4).

\end{document}